\documentclass[12pt,a4paper]{article}
\usepackage{theorem}
\newtheorem{lemma}{Lemma}

\newtheorem{theorem}[lemma]{Theorem}

{\theorembodyfont{\upshape}}
{\theorembodyfont{\upshape}}
{\theorembodyfont{\upshape}}
{\theorembodyfont{\upshape}}
\newcommand{\R}{{\bf R}}
\newcommand{\C}{{\bf C}}
\newcommand{\rme}{{\rm e}}
\newcommand{\rmd}{{\rm d}}
\newcommand{\sig}{\sigma}

\newcommand{\bet}{\beta}
\newcommand{\gam}{\gamma}

\newcommand{\del}{\delta}
\newcommand{\eps}{\varepsilon}

\newcommand{\lap}{{\Delta}}

\newcommand{\Spec}{{\rm Spec}}

\newcommand{\norm}{\Vert}
\renewcommand{\Re}{{\rm Re}}
\renewcommand{\Im}{{\rm Im}}

\newcommand{\Proof}{\underbar{Proof}{\hskip 0.1in}}

\newcommand{\Schrodinger}{Schr\"odinger }
\newcommand{\pr}{\prime}
\title{SEMI-CLASSICAL STATES \break
FOR NON-SELF-ADJOINT\break
SCHR\"ODINGER OPERATORS}
\author{E.B. Davies}
\date{March 1998}
\begin{document}
\maketitle
%%%%%%%%%%%%%%%%%%%%%
\begin{abstract}
We prove that the spectrum of certain  
non-self-adjoint \Schrodinger operators is unstable in the 
semi-classical limit $h\to 0$. Similar 
results hold for a fixed operator in the high energy limit. 
The method involves the construction of approximate semi-classical 
modes of the operator by the JWKB method for energies far 
from the spectrum. 
%%%%%%%%%%%%%%%%%%
\vskip 0.1in
keywords: semiclassical limit, JWKB method, 
non-self-adjoint operators, spectrum, pseudospectrum, resolvent norm,
spectral instability.
\end{abstract}
%%%%%%%%%%%%%%%%%%%%%%%%
\section{Introduction}
\par
It is well known that the complex spectrum of many
non-self-adjoint differential operators is highly unstable 
under small perturbations \cite{D1, D2, RT, T}; this 
has been investigated in detail for the Rayleigh equation 
in hydrodynamics in \cite[Ch 4]{KS}. One way 
of exploring this 
fact is by defining the pseudospectrum of such an 
operator $H$ by
\[
\Spec_{\eps}(H):=\Spec(H)\cup \lbrace z: \norm 
R(z)\norm>\eps^{-1}\rbrace
\]
where $\eps >0$ and $R(z)$ is the resolvent of $H$.
It is known that $\Spec_{\eps}(H)$ contains the 
$\eps$-neighbourhood of the spectrum, and that it is 
contained in the $\eps$-neighbourhood of the numerical 
range of $H$. Theorem 1 states that for \Schrodinger 
operators with complex potentials 
$\Spec_{\eps}(H)$ expands to fill a region $U$ of the complex 
plane much larger than the spectrum in the semi-classical 
limit $h\to 0$. More 
precisely we obtain an explicit lower bound on $\norm 
R(z)\norm $ which increases rapidly as $h\to 0$ for 
all $z$ in the region $U$ defined by
\[
U:=\lbrace z=\eta^{2}+V(a):\eta\in\R\backslash \lbrace 
0\rbrace  \mbox{   and   } 
\Im(V^{\pr}(a))\not= 0\rbrace.
\]
For $h=1$ we deduce that for large enough $z$ within a 
suitable region the resolvent norms of \Schrodinger 
operators with complex potentials become very large even 
though $z$ may be far 
from the spectrum. This had apparently not been noticed in 
other spectral investigations of these operators \cite{AMT, 
DR, EE, LM}, 
apart from \cite{D2}, whose results are greatly extended and 
improved in Theorem 2 below.

Our results have a positive aspect. Our proofs use a JWKB 
analysis to construct a continuous family of approximate 
eigenstates, which we call semi-classical modes, for 
the operators in question. These modes have complex energies 
far from the spectrum, but could be used to investigate the 
time evolution of fairly general initial states by 
expanding these in terms of the modes.

%%%%%%%%%%%%%%%%%%%%%%%%%%%%%%%%%%%%%
\section{The estimates}
For reasons which will become clear in the next section, we 
consider operators somewhat more general than those 
described in the last section. We assume that $H$ acts in 
$L^{2}(\R)$ according to the formula 
\[
Hf(x):=-h^{2}\frac{\rmd^{2} f}{\rmd x^{2}} +V_{h}(x)f(x)
\]
Where $V_{h}$ are smooth potentials for all small enough $h>0$ 
which depend continuously on $h$, together with 
their derivatives of all orders. In the applications in the 
next section $V_{h}$ 
has an expansion involving fractional powers of $h$, 
but this is invisible since we treat $V_{h}$ as it stands, 
and only need an asymptotic expansion involving integer 
powers of $h$; this simplification is essential to our 
solution of the problem. Technically 
we assume that $H$ is 
some closed extension of the operator initially defined on 
$C_{c}^{\infty}(\R)$. We assume throughout the 
section that 
\[
z:=\eta^{2}+V_{h}(a)
\]
where $\eta\in\R\backslash \lbrace 0\rbrace$ and 
$\Im (V_{0}^{\pr}(a))\not= 0$. The $h$-dependence of $z$ as 
defined above 
may be eliminated by a uniformity argument spelled out in 
the next section, but we do not focus on this issue here. 
Our goal is to prove the upper bound
\[
\norm H\tilde f -z\tilde f\norm/\norm \tilde f \norm =O(h^{n})
\]
as $h\to 0$, for all $n>0$, where $\tilde f\in 
C_{c}^{\infty}(\R)$ depends upon $h$ and $n$. This immediately 
implies that $\norm (H-zI)^{-1}\norm$ diverges as $h\to 0$ 
faster than any negative power of $h$. Although the above 
equation may be interpreted as stating that such $z$ are 
approximate eigenvalues for small enough $h>0$, it does 
not follow that they are 
close to true eigenvalues, and indeed the examples studied 
in \cite{D1, D2,RT,T} show that there is a strong distinction 
between the spectrum and pseudospectrum.
\par
%%%%%%%%%%%%%%%%%%%%%%%%%%%%
\begin{theorem}
There exists $\del >0$ and for each $n>0$ a positive 
constant $c_{n}$ and an $h$-dependent function $\tilde f\in 
C_{c}^{\infty}(\R)$ such that if $0<h<\del$ then
\[
\norm H\tilde f -z\tilde f\norm/\norm \tilde f \norm \leq c_{n}h^{n}.
\]
\end{theorem}
\Proof The proof is a direct construction. We put $\tilde 
f(a+s):=\xi(s) f(s)$ for all $s\in \R$ 
where $\xi\in C_{c}^{\infty}(\R)$ satisfies $\xi(s)=1$ if 
$|s|<\del/2$ and $\xi(s)=0$ if $|s|>\del$, and $\del>0$ is 
determined below. We take $f$ to be the smooth but not 
square integrable JWKB function 
$f:=\exp(-\psi)$ where
\[
\psi(s):=\sum_{m=-1}^{n}h^{m}\psi_{m}(s)
\]
and $\psi_{m}$ are defined below. A direct computation 
shows that 
\[
Hf-zf=\left(\sum_{m=0}^{2n+2}h^{m}\phi_{m}\right) f
\]
where $\phi_{m}$ are given by the formulae
\begin{eqnarray*}
\phi_{0}&:=&-(\psi_{-1}^{\pr})^{2}+V_{h}-z\\
\phi_{1}&:=&\psi_{-1}^{\pr\pr}-2\psi_{-1}^{\pr}\psi_{0}^{\pr}\\
\phi_{2}&:=&\psi_{0}^{\pr\pr}-2\psi_{-1}^{\pr}\psi_{1}^{\pr}
-(\psi_{0}^{\pr})^{2}\\
\phi_{3}&:=&\psi_{1}^{\pr\pr}-2\psi_{-1}^{\pr}\psi_{2}^{\pr}
-2\psi_{0}^{\pr}\psi_{1}^{\pr}\\
&\ldots&\\
\phi_{2n+2}&:=&-(\psi_{n}^{\pr})^{2}.\\
\end{eqnarray*}
By setting $\phi_{m}=0$ for $0\leq m\leq n+1$, we obtain a 
series of equations which enable us to determine all 
$\psi_{m}$, provided $\del>0$ is small enough.
\par
The key equation is the complex eikonal equation
\[
\psi_{-1}^{\pr}(s)^{2}=V_{h}(a+s)-V_{h}(a)-\eta^{2}
\]
whose solution is
\begin{eqnarray*}
\psi_{-1}(s)&:=&\int_{0}^{s}
\left( V_{h}(a+t)-V_{h}(a)-\eta^{2}\right)^{1/2} \rmd t\\
&=& \int_{0}^{s}i\eta 
\left(1-\frac{ V_{h}(a+t)-V_{h}(a)}{\eta^{2}}\right)^{1/2} \rmd 
t\\
&=& i\eta s-\frac{iV_{h}^{\pr}(a)s^{2}}{4\eta} +O(s^{3}).
\end{eqnarray*}
We have assumed that $\Im(V_{0}^{\pr}(a))\not= 0$ and this 
implies that $\Im(V_{h}^{\pr}(a))\not= 0$ for all small enough 
$h>0$. We choose $\eta$ to be of the same sign as 
$\Im(V_{h}^{\pr}(a))$ so that for a suitable 
constant $\gam>0$ we have 
\[
\gam s^{2} \leq \Re (\psi_{-1}(s)) \leq 3\gam s^{2}
\]
for all small enough $s$ and $h$. We also assume that $s$ 
and $h$ are small enough for 
$\rho:=(2\psi_{-1}^{\pr})^{-1}$ to satisfy a bound of the 
form $|\rho(s)|\leq \bet$.
\par
We now force $\phi_{m}=0$ for $0\leq m\leq n+1$ by putting
\begin{eqnarray*}
\psi_{0}^{\pr}&=& \rho \psi_{-1}^{\pr\pr}\\
\psi_{1}^{\pr}&=& \rho 
(\psi_{0}^{\pr\pr}-(\psi_{0}^{\pr})^{2})\\
\psi_{2}^{\pr}&=& \rho 
(\psi_{1}^{\pr\pr}-2\psi_{0}^{\pr}\psi_{1}^{\pr})\\
&\mbox{etc.}&
\end{eqnarray*}
We determine the functions uniquely by also imposing 
$\psi_{m}(0)=0$ for all $0\leq m\leq n$. 
Each of the functions is bounded provided $s$ and $h$ are 
small enough, and the same is true of the remaining 
functions $\phi_{m}$. Specifically we assume that for 
some $\del >0$ and constants $c_{m},c_{m}^{\pr}$ we have
\[
|\psi_{m}(s)|\leq c_{m}
\]
for $0\leq m\leq n$, and
\[
|\phi_{m}(s)|\leq c_{m}^{\pr}
\]
for $n+2\leq m\leq 2n+2$, provided $|s|\leq \del$ and $0<h<\del^{2}$.
\par
In the following calculations $a_{i}$ denote various 
positive constants, independent of $h$ and $s$. We have
\begin{eqnarray*}
\norm \tilde f \norm_{2}^{2}&\geq &
\int_{-\del/2}^{\del/2}|f(s)|^{2}\rmd s\\
&\geq & 
\int_{-\del/2}^{\del/2}\rme^{-3\gam s^{2}h^{-1}-a_{1}}\rmd s\\
&=&\int_{-\del h^{-1/2}/2}^{\del h^{-1/2}/2}
\rme^{-3\gam t^{2}-a_{1}} h^{1/2}\rmd t\\
&\geq & \int_{-1/2}^{1/2}
\rme^{-3\gam t^{2}-a_{1}}h^{1/2} \rmd t\\
&=& a_{2}h^{1/2}.
\end{eqnarray*}
We also have
\begin{eqnarray*}
\norm H\tilde f -z\tilde f \norm_{2} &=&\norm -h^{2}f 
\xi^{\pr\pr} -2h^{2} f^{\pr}\xi^{\pr} +\xi(Hf-zf)\norm_{2}\\
&\leq & h^{2}\norm f \xi^{\pr\pr}\norm_{2} +
2h^{2} \norm f^{\pr}\xi^{\pr}\norm_{2} +
\sum_{m=n+2}^{2n+2}h^{m}\norm \xi \phi_{m} f \norm_{2}
\end{eqnarray*}
and need to estimate each of the norms. Since $\xi^{\pr}$ 
has support in $\lbrace s:\del/2 \leq |s|\leq \del 
\rbrace$, we have
\begin{eqnarray*}
\norm \xi^{\pr\pr}f\norm_{2}^{2}&\leq&
a_{3}\int_{\del/2 \leq |s|\leq \del} \rme^{-\gam 
s^{2}h^{-1}+a_{4}}\rmd s\\
&\leq& a_{5}\rme^{-\gam\del^{2}/4h}.
\end{eqnarray*}
In other words  $\norm \xi^{\pr\pr}f\norm_{2}$ decreases at an 
exponential rate as $h\to 0$. A similar argument applies to 
$\norm \xi^{\pr}f^{\pr}\norm_{2}$. Since $\phi_{m}$ is bounded 
on $\lbrace s:|s|\leq \del\rbrace$, uniformly for $|h|\leq 
\del^{2}$, we see that 
\begin{eqnarray*}
\norm \xi \phi_{m} f \norm_{2}^{2}&\leq & a_{6}
\int_{-\del}^{\del}|f(s)|^{2}\rmd s\\
&\leq & 
\int_{-\del}^{\del}\rme^{-\gam s^{2}h^{-1}+a_{7}}\rmd s\\
&\leq & a_{8}h^{1/2}
\end{eqnarray*}
by an argument similar to that used for $\tilde f$ above.
Putting the various inequalities together we obtain
the statement of the theorem.
\par
%%%%%%%%%%%%%%%%%%%%%%%%%%%
\section{High Energy Spectrum}
\par
By a change of scale our theorem can be applied to prove 
the instability of the high energy spectrum of a 
non-self-adjoint \Schrodinger operator with complex 
potential. The results in this section extend those of 
\cite{D2} both by providing greater insight into the 
machanism involved and by obtaining much stronger estimates.  
Adopting quantum mechanical notation we assume 
that $h=1$ and that the operator $H$ acting in $L^{2}(\R)$ 
is given by 
\[
H:=P^{2}+\sum_{m=1}^{n}c_{m}Q^{m}
\]
where $n$ is even and the constant $c_{n}$ has positive 
real and imaginary parts. 
\par
%%%%%%%%%%%%%%%%%%%%%%%%%%%%
\begin{theorem}
If $z\in \C$ satisfies $0<\arg(z)<\arg(c_{n})$ and $\sig >0$ then
\[
\norm (H-\sig zI)^{-1}\norm
\]
diverges to infinity faster than any power of $\sig$ as 
$\sig\to +\infty$.
\end{theorem}
\par
\Proof If $u>0$  then the operator $H$ is unitarily 
equivalent to the operator
\[
H_{1}:=u^{-2}P^{2}+\sum_{m=1}^{n}c_{m}u^{m}Q^{m}
\]
Putting $u:=\sig^{1/n}$ we obtain
\[
\norm (H-\sig zI)^{-1}\norm =\sig^{-1} \norm (H_{2}- zI)^{-1}\norm
\]
where
\[
H_{2}:=\sig^{-1}H_{1}=u^{-2-n}P^{2}+\sum_{m=1}^{n}c_{m}u^{m-n}Q^{m}
\]
Putting $h:=u^{-(n+2)/2}$ we have
\[
H_{2}:=h^{2}P^{2}+\sum_{m=1}^{n}c_{m}h^{2(n-m)/(n+2)}Q^{m}
=h^{2}P^{2}+V_{h}(Q).
\]
This is precisely the form of operator to which Theorem 1 
applies. We have $V_{0}^{\pr}(a)=c_{n}a^{n}$ so 
$\Im(c_{n})>0$ implies that the conditions of Theorem 1 are 
satisfied for any $z$ in the sector
\[
U:=\lbrace z: 0<\arg (z) < \arg (c_{n})\rbrace.
\]
This completes the proof, except for a technical point 
which we now address. In Theorem 1 we assumed that 
$z=\eta^{2}+V_{h}(a)$ where $\Im(V_{0}^{\pr}(a))\not= 0$, 
so $z$ is apparently dependent on $h$, with the limit 
$z_{0}:=\eta^{2}+V_{0}(a)$ as $h\to 0$. We rectify this 
problem by fixing $z\in U$ and making $a$ and $\eta$ depend upon 
$h$ in
such a way that
\[
z=\eta_{h}^{2}+V_{h}(a_{h})
\]
where $\eta_{h}\to \eta$ and $a_{h}\to a$ as $h\to 0$. We 
now have to check that all the estimates of Section 2 are 
locally uniform with respect to $\eta$ and $a$, so that the 
result we claim does indeed follow.
\par
The method of this paper can be extended to treat certain
rotationally invariant problems in higher space dimensions. 
The condition $-2\leq p(1)$ in the next theorem is included 
because it is relevant to the existence of a closed 
extension of the operator, by virtue of an application of 
the theory of sectorial forms.
%%%%%%%%%%%%%%%%%%%%%%%%
\begin{theorem}
Let the operator $H$ acting in $L^{2}(\R^{N})$ be some 
closed extension of the operator given by
\[
Hf(x):=-\lap f(x) +\sum_{m=1}^{n}c_{m}|x|^{p(m)}f(x)
\]
for all $f$ in the initial domain 
$C_{c}^{\infty}(\R^{N}\backslash \lbrace 0\rbrace)$, where $c_{n}$ has 
positive real and imaginary parts, $p(n)>0$ and
\[
-2\leq p(1)<p(2)<\ldots<p(n).
\]
If $z\in \C$ satisfies $0<\arg(z)<\arg(c_{n})$ and $\sig >0$ then
\[
\norm (H-\sig zI)^{-1}\norm
\]
diverges to infinity faster than any power of $\sig$ as 
$\sig\to +\infty$.
\end{theorem}
\par
\Proof The difference from Theorem 2 is that 
after restricting to the usual angular momentum sectors the 
operators act in $L^{2}(0,\infty)$ and include angular 
momentum terms in the potential. However, it may be seen 
that the analysis of Theorem 2 can be extended to operators of the form
\[
H:=P^{2}+\sum_{m=1}^{n}c_{m}Q^{p(m)}
\] 
so the incorporation of the 
angular momentum terms causes no difficulties.
\par
\underbar {Note} Since the supports of the test functions used in the proof of 
the theorem are compact and move to infinity, Theorem 3 
remains valid if we add a non-central potential to $H$, 
provided that potential decreases at infinity faster than 
any negative power of $|x|$. Weaker versions of Theorems 2 and 3 hold if 
one adjoins a potential which decreases more slowly at 
infinity. 
%%%%%%%%%%%%%%%%%%%%%%%%%
\vskip 0.3in
{\bf Acknowledgments } I should like to thank M Kelbert and 
Y  Safarov for helpful comments.
%\vskip 0.3in
\newpage
%%%%%%%%%%%%%%%%%%%%%%%

\par
%%%%%%%%%%%%%%%%%%%%%%%%%%%%%%%%%%%%%%%%
\vskip 0.3in
Department of Mathematics \newline
King's College \newline
Strand \newline
London WC2R 2LS \newline
England \newline
e-mail: E.Brian.Davies@kcl.ac.uk
\end{document}